\magnification=1100

%

%
%
\font\seventeenbf=cmbx10 at 17.28pt      
\font\seventeenbfit=cmmib10 at 17.28pt  
\font\seventeenbsy=cmbsy10 at 17.28pt
%
%
\font\fourteenbf=cmbx10 at 14.4pt
\font\fourteenbfit=cmmib10 at 14.4pt
\font\fourteenbsy=cmbsy10 at 14.4pt
%
%
\font\twelvebf=cmbx10 at 12pt
\font\twelvebfit=cmmib10 at 12pt
\font\twelvebsy=cmbsy10 at 12pt 
%
%
%
%
\font\tenbsy=cmbsy10

\font\eightbsy=cmbsy8
\font\sevenbsy=cmbsy7
\font\sixbsy=cmbsy6
\font\fivebsy=cmbsy5
\newfam\bsyfam
  \textfont\bsyfam=\tenbsy
  \scriptfont\bsyfam=\sevenbsy
  \scriptscriptfont\bsyfam=\fivebsy
 
%
%
\font\twelvemsa=msam10 at 12pt
\font\fourteenmsa=msam10 at 14.4pt
\font\seventeenmsa=msam10 at 17.28pt
\font\tenmsa=msam10

\font\sevenmsa=msam7
\font\fivemsa=msam5
\newfam\msafam
  \textfont\msafam=\tenmsa
  \scriptfont\msafam=\sevenmsa
  \scriptscriptfont\msafam=\fivemsa
 
%
%
\font\twelvemsb=msbm10 at 12pt
\font\fourteenmsb=msbm10 at 14.4pt
\font\seventeenmsb=msbm10 at 17.28pt
\font\tenmsb=msbm10
\font\ninemsb=msbm9
\font\sevenmsb=msbm7
\font\fivemsb=msbm5
\newfam\msbfam
  \textfont\msbfam=\tenmsb
  \scriptfont\msbfam=\sevenmsb
  \scriptscriptfont\msbfam=\fivemsb

\def\Bbb{\fam\msbfam\tenmsb}
%
%
\font\twelveCalbf=eusb10 at 12pt
\font\fourteenCalbf=eusb10 at 14.4pt
\font\seventeenCalbf=eusb10 at 17.28pt
\font\nineCal=eusm9
\font\tenCal=eusm10

\font\sevenCal=eusm7
\font\fiveCal=eusm5
\newfam\Calfam
  \textfont\Calfam=\tenCal
  \scriptfont\Calfam=\sevenCal
  \scriptscriptfont\Calfam=\fiveCal
\def\Cal{\fam\Calfam\tenCal}
%
%
\font\twelveeuf=eufm10 at 12pt
\font\fourteeneuf=eufm10 at 14.4pt
\font\seventeeneuf=eufm10 at 17.28pt
\font\teneuf=eufm10
\font\seveneuf=eufm7
\font\fiveeuf=eufm5
\newfam\euffam
  \textfont\euffam=\teneuf
  \scriptfont\euffam=\seveneuf
  \scriptscriptfont\euffam=\fiveeuf

%
%
\font\tenbfit=cmmib10
\font\eightbfit=cmmib8
\font\sevenbfit=cmmib7
\font\sixbfit=cmmib6
\font\fivebfit=cmmib5
\newfam\bfitfam
  \textfont\bfitfam=\tenbfit
  \scriptfont\bfitfam=\sevenbfit
  \scriptscriptfont\bfitfam=\fivebfit

%
%
\font \seventeensans          = cmss10 at 17.28pt

\font \tensans                = cmss10
\font \fivesans               = cmss10 at 5pt
\font \sixsans                = cmss10 at 6pt
\font \sevensans              = cmss10 at 7pt
\font \ninesans               = cmss10 at 9pt
\newfam\sansfam
        \textfont\sansfam=\tensans
        \scriptfont\sansfam=\sevensans
        \scriptscriptfont\sansfam=\fivesans
\def\sans{\fam\sansfam\tensans}
%
%
\font\sixrm=cmr6
\font\sixi=cmmi6
\font\sixsy=cmsy6
\font\sixbf=cmbx6
%
%

\font\eightbf=cmbx8

%
%
\font \ninerm                 = cmr9
\font \ninei                  = cmmi9
\font \ninebf                 = cmbx9
\font \ninesy                 = cmsy9
\font \nineit                 = cmti9
\font \ninett                 = cmtt9
\font \ninesl                 = cmsl9
%
%

\font \twelvebf               =cmbx10 at 12pt
\font \twelvebfit             =cmmib10 at 12pt

\font \twelvebsy              =cmbsy10 at 12pt 

\def\smallfonts{ 
 \def\rm{\fam0\ninerm}%
                \textfont0=\ninerm   \scriptfont0=\sixrm    \scriptscriptfont0=\fiverm
                \textfont1=\ninei    \scriptfont1=\sixi     \scriptscriptfont1=\fivei
                \textfont2=\ninesy   \scriptfont2=\sixsy    \scriptscriptfont2=\fivesy
 \def\it{\fam\itfam\nineit}
                \textfont\itfam=\nineit
 \def\sl{\fam\slfam\ninesl}
   \textfont\slfam=\ninesl
 \def\bf{\fam\bffam\ninebf}
   \textfont\bffam=\ninebf \scriptfont\bffam=\sixbf
   \scriptscriptfont\bffam=\fivebf
 \def\sans{\fam\sansfam\ninesans}%
    \textfont\sansfam=\ninesans \scriptfont\sansfam=\sixsans
    \scriptscriptfont\sansfam=\fivesans
 \def\tt{\fam\ttfam\ninett}
    \textfont\ttfam=\ninett
 \def\Bbb{\fam\msbfam\ninemsb}
  \textfont\msbfam=\ninemsb
 \def\Cal{\fam\Calfam\nineCal}
  \textfont\Calfam=\nineCal
\rm}%
\def\small{
  \vskip6pt\bgroup\smallfonts 
   \normalbaselineskip=10pt\parindent=0pt\parskip=0pt%
\noindent
}
\def\endsmall{\vskip6pt\egroup%
}


\def\seventeenpointbf{%
 \textfont0=\seventeenbf  \scriptfont0=\twelvebf  \scriptscriptfont0=\eightbf
 \textfont1=\seventeenbfit\scriptfont1=\twelvebfit\scriptscriptfont1=\eightbfit
 \textfont2=\seventeenbsy \scriptfont2=\twelvebsy \scriptscriptfont2=\eightbsy
 \textfont\msafam=\seventeenmsa    
 \textfont\msbfam=\seventeenmsb    
 \textfont\sansfam=\seventeensans  
 \textfont\euffam=\seventeeneuf    
 \textfont\Calfam=\seventeenCalbf
 \seventeenbf
 \normalbaselineskip=20.736pt\normalbaselines}


\def\fourteenpointbf{%
 \textfont0=\fourteenbf   \scriptfont0=\tenbf   \scriptscriptfont0=\sevenbf
 \textfont1=\fourteenbfit \scriptfont1=\tenbfit \scriptscriptfont1=\sevenbfit
 \textfont2=\fourteenbsy  \scriptfont2=\tenbsy  \scriptscriptfont2=\sevenbsy
 \textfont\msafam=\fourteenmsa    
 \textfont\msbfam=\fourteenmsb    
 \textfont\euffam=\fourteeneuf    
 \textfont\Calfam=\fourteenCalbf
 \fourteenbf
 \normalbaselineskip=17.28pt\normalbaselines}


\def\twelvepointbf{%
 \textfont0=\twelvebf   \scriptfont0=\eightbf   \scriptscriptfont0=\sixbf
 \textfont1=\twelvebfit \scriptfont1=\eightbfit \scriptscriptfont1=\sixbfit
 \textfont2=\twelvebsy  \scriptfont2=\eightbsy  \scriptscriptfont2=\sixbsy
 \textfont\msafam=\twelvemsa    
 \textfont\msbfam=\twelvemsb    
 \textfont\euffam=\twelveeuf    
  \textfont\Calfam=\twelveCalbf
 \twelvebf
 \normalbaselineskip=14.4pt\normalbaselines}


\def\tenpointbf{%
 \textfont0=\tenbf   \scriptfont0=\sevenbf   \scriptscriptfont0=\fivebf
 \textfont1=\tenbfit \scriptfont1=\sevenbfit \scriptscriptfont1=\fivebfit
 \textfont2=\tenbsy  \scriptfont2=\sevenbsy  \scriptscriptfont2=\fivebsy
 \textfont\msbfam=\tenmsb       
 \tenbf}


\def\titlefont{\seventeenpointbf}
\def\sectionfont{\twelvepointbf}
\def\proclaimfont{\tenpointbf}
\def\chapterfont{\fourteenpointbf}




\def\AA{{\Cal A}}  
  
   
\def\FF{{\Cal F}}   
\def\GG{{\Cal G}}


\baselineskip=12pt     
\parskip=5.5pt         

\newskip\baseskip
\baseskip=17.5pt plus 2pt minus 1pt  

\newskip\skchpre
\newskip\skchpost
\newskip\sksecpre
\newskip\sksecpost

\skchpre=2\baseskip
\skchpost=5\baseskip
\sksecpre=\baseskip
\sksecpost=0.7\baseskip

\def\postchapskip{\nobreak\vskip\skchpost\nobreak}

\def\presecskip{\smallbreak\vskip\sksecpre\nobreak}
\def\postsecskip{\nobreak\vskip\sksecpost\nobreak}

\def\chapter#1{\vfill\supereject
   \noindent{{\chapterfont Chapter
        \hglue-2pt\newchapternumber\hglue-4pt.~~#1} } 
                        \postchapskip}
                          

\newcount\labelREG
    \def\newlabelnumber{%
    \global\advance\labelREG by1
    \the\labelREG%
       }
    \def\labelnumber{\the\labelREG}

\newcount\chapREG  
     \def\newchapternumber{ 
     \global\advance\chapREG by1
     \global\secREG=0   
     \global\labelREG=0 
     \uppercase\expandafter{\romannumeral\the\chapREG}
     } 
     \def\chapternumber{\uppercase\expandafter{\romannumeral\the\chapREG}}

\def\section#1{
        \presecskip
   \noindent%
   {\S \sectionfont\newsectionnumber.\quad #1}
   \postsecskip} 

     \newcount\secREG  
     \def\newsectionnumber{ 
     \global\advance\secREG by1
     \global\labelREG=0
     \the\secREG
     } 
    \def\sectionnumber{\the\secREG}
    \def\setsection#1{\global\secREG=#1}

 
\def\proclaim#1 
 {\medbreak 
  \smallskip 
  \noindent 
   {\proclaimfont \sectionnumber.\newlabelnumber.\hglue3pt #1.} 
   \begingroup\sl\penalty80 
} 
\def\endproclaim{\endgroup\smallskip} 

\def\nnproclaim#1 
  {\medbreak 
   \smallskip 
   \noindent 
   {\bf#1.}\   
   \begingroup\sl\penalty80
  } 
 \def\endnnproclaim{\endgroup\smallskip}                      

\def\pproclaim#1{
  \proclaim#1
  \rm
}
\def\endpproclaim{
   \endproclaim
}

\def\proof{
  \noindent
  {\bf Proof:}
}

\def\endproof{
  Q.E.D.\medskip
}



\newcount\itemregister  
 
 
\def\itemnumber{ 
      {\global\advance\itemregister by1} 
      \the\itemregister}                 
 

\def\mapright#1{\mathop{\vbox{\ialign{
                                ##\crcr
    ${\scriptstyle\hfil\;\;#1\;\;\hfil}$\crcr
 \noalign{\kern2pt\nointerlineskip}
    \rightarrowfill\crcr}}\;}}

\def\mapleft#1{\mathop{\vbox{\ialign{
                                ##\crcr
    ${\scriptstyle\hfil\;\;#1\;\;\hfil}$\crcr
    \noalign{\kern2pt\nointerlineskip}
    \leftarrowfill\crcr}}\;}}



\def\perp{^\bot} 


\def\ie{{\it i.e.},\ }

\def\comp{\raise1pt\hbox{{$\scriptscriptstyle\circ$}}}

\def\dbar{{\overline{\partial}}}

\def\deg{\mathop{\rm deg\, }\nolimits} 
\def\del{\partial}

\def\dim{\mathop{\rm dim\, }\nolimits} 
\def\End{\mathop{\rm End}\nolimits} 
\def\exp{\mathop{\rm exp}\nolimits} 
 
\def\genus{\mathop{{\rm genus}}\nolimits}

\def\hdg{{\rm Hdg}}
 
\def\ii{{\rm i}}

\def\Im{\mathop{\rm Im}\nolimits}

\def\Ker{\mathop{\rm Ker}\nolimits}
\def\ker{\mathop{\rm Ker}\nolimits}

\def\rank{\mathop{\rm rank}\nolimits}


\def\SL{\mathop{\rm SL}\nolimits}


\def\Tr{{}^T\kern-0.9pt} 
\def\tr{{}^t\kern-0.9pt} 
\def\trace{\mathop{\rm Tr}\nolimits}

\def\ve{{\scriptscriptstyle \vee}}

%
%
%

\line{\hfil\titlefont  Arakelov-type inequalities\hfil} 
\line{\hfil\titlefont for Hodge bundles\hfil}
\vskip 3em
\line{\hfil Chris Peters\hfil}
\line{\hfil\tt Department of Mathematics\hfil}
\line{\hfil\tt University of Grenoble I\hfil}
\line{\hfil\tt Saint-Martin d'H{\`e}res, France\hfil}
\vskip 3em
\line{\hfill \bf June 20 2000\hfill }
\vskip 2em 
\line{\hfill Pr{\'e}publication de l'Institut Fourier n$^o$ 511  (2000)\hfill}
\line{\hfill \sl
http://www-fourier.ujf-grenoble.fr/prepublications.html\hfill} 
\footnote{}{\small  {\bf Keywords:} Variations of Hodge structure,
Higgs bundles, Higgs field, period map, Hodge bundles.\break
  \endsmall}
\setsection{-1}
\section{Introduction}

The inequalities from the title refer back to Arakelov's article
[Arakelov]. The main result of that paper is:
\nnproclaim{Theorem} Fix a complete curve $C$ of genus  $>1$ and a 
finite set $S$ of points on $C$. There are at most finitely
many non-isotrivial families of curves of given genus over $C$ that 
are  smooth over $C\setminus S$. 
\endnnproclaim

The proof consists of two parts. First one proves that there are only
finitely many such families (this is a {\it boundedness} statement) by
bounding the degree $d$ of the relative canonical bundle in terms of
the genus $p$ of $C$, the genus $g$ of the fiber and the cardinality
of the set $S$~:
$$
0\le d \le (2p-2+\# S){g\over 2}.
$$
The second part consists of establishing {\it rigidity} for a
non-isotrivial  family. It follows upon identifying the deformation
space of the family  with  the $H^1$ of the inverse of the relative
canonical bundle, which  is shown to be ample. Kodaira vanishing then
completes the proof.

This approach can be carried out for other situations 
as well. In fact [Faltings] deals with the case of abelian varieties and 
shows that boundedness always holds and that for rigidity one has
to impose further conditions besides non-isotriviality. 
Subsequently the rigidity statement  has been  generalized in [Peters90] 
and  using his result, the case of K3-surfaces, resp. Abelian varieties could
be treated completely by Saito and Zucker in [Saito-Zucker],
resp. by Saito in [Saito].  

The boundedness statement is an inequality for the
degree $d$ of the direct image of the relative canonical bundle, \ie the
(canonical extension) of the Hodge bundle and one can ask for bounds
for the degrees of the other Hodge bundles. In fact, the main result
of this note gives such a bound for the Hodge components of complex
variations of Hodge structures in terms of ranks of iterates of the
Higgs field (the linear map between Hodge components induced by the
Gauss-Manin connection). By way of an example, we have~:

\nnproclaim{Proposition} Let $V=\bigoplus_{p=0}^w V^{p,w-p}$ be a real
weight $w$ variation of Hodge structures over a punctured curve
$C\setminus S$ with unipotent local monodromy-operators. Let
$\sigma^p:V^{p,w-p}\to V^{p-1,w-p+1}$ the $k$-th component of the
Higgs field and $\sigma^k=\sigma_{w-k+1}\comp\cdots\comp
\sigma_w:V^{w,0}\to V^{w-k,k}$ the $k$-th iterate. Then
$$
0\le \deg V^{w,0} \le (2p-2+\# S)\left(\sum_{r=1}^w {r\over 2}
\left(\rank\sigma^{r-1}-\rank \sigma^r\right)\right).
$$
\endnnproclaim

Some time ago I sketched a proof of a similar, but weaker inequality
in [Peters86], but the details of this proof never appeared.  Deligne
found an amplification of my argument when the base is a {\it compact}
curve leading to optimal bounds for {\it complex} variations of Hodge
structures (letter to the author 18/2/1986).  The principal goal of
this note is to give a complete proof of the refined inequalities (for
a complex variation with quasi-unipotent local   monodromy-operators over a
quasi-projective smooth curve) based on this letter in the
light of later developments which I sketch below. 

The reason for writing up this note stems from a recent revival of
interest in this circle of ideas: on the one hand Parshin posed me some
questions related to this. On the other hand, Jost and Zuo sent me a
preprint [JostZuo] containing similar (but weaker) bounds obtained by
essentially the same method.

Continuing with the historical development, Deligne's letter and
subsequent correspondence between Deligne and Beilinson together with
basic ideas and results of Hitchin paved the way for the theory of
Higgs bundles, the proper framework for such questions.  See
[Simpson92], [Simpson94] and [Simpson95] for a further discussion of
these matters. 

In relation with Simpson's work, I should remark that the boundedness
result from [Simpson94] (Corollary 3.4 with $\mu=0$ and $P=0$)
immediately implies that the Chern numbers of Hodge bundles underlying
a complex variation of given type over a compact projective variety
can only assume finitely many values.  It follows that there are {\sl
a priori} bounds on these Chern numbers. If the base is a (not
necessarily compact) curve, Simpson shows in [Simpson90] that the same
methods give boundedness for the canonical extensions of the Hodge
bundles. Explicit bounds were not given however.  Conversely, knowing that
bounds on the degrees of the Hodge bundles exist can be used to
simplify the rather technical proof for boundedness, as outlined in \S 4.

As to the further contents of this note, in \S1 I rephrase known bounds on
the curvature of Hodge bundles (with respect to the Hodge metric) in a
way that shows how to adapt these in the case of non-compact algebraic
base manifolds. The main argument is given in \S2 while in \S3 I
explain what one has to change  in the non-compact case. 

In closing, I want to mention Eyssidieux work which treats a
generalization of the true Arakelov inequality (weight one Hodge
structures on a compact curve) to higher dimensional base
manifolds. The inequalities concern variations of Hodge structure over
a compact K{\"a}hler manifold $M$ such that the period map is
generically finite onto its image. See [Eyss].

\section{Curvature bounds for Hodge bundles}

Let me recall that a {\it complex variation of Hodge bundles of weight
$w$} on a complex manifold $M$ consists of a complex local system $V$
on $M$ with a direct sum decomposition into complex subbundles 
$$
V=\bigoplus_{p+q=w} V^{p,q}
$$
with the property that the canonical flat connection $\nabla$
satisfies the transversality condition
$$
\nabla : V^{p,q} \to \AA^{1,0}(V^{p-1,q+1})\oplus
\AA^1(V^{p,q})\oplus\AA^{0,1}(V^{p+1.q-1}). \leqno{\rm(*)}
$$
The local system $V$ defines a holomorphic vector bundle denoted by
the same symbol. The bundles $V^{p,q}$ are not necessarily
holomorphic, but the transversality condition implies that the filtration defined by
$$
F^pV:=\bigoplus _{r\ge p} V^{r,w-r}
$$
is holomorphic. We have a $C^\infty$-isomorphism between the {\it Hodge bundle}
$$
V_\hdg:= \bigoplus_p F^{p+1}V/F^pV 
$$
and the bundle $V$ which we shall use to transport the connection to the 
former. 

Of course, the usual (real) variations are examples of complex
variations. These satisfy the additional {\it reality constraint}
$V^{p,q}=\bar V^{q,p}$.  Conversely, a complex variation $V$ together
with its complex conjugate, defines a real variation of Hodge
structures on $V\oplus \bar V$ in the obvious manner.  In passing, we
observe however that we can have complex variations of any given pure
type.  If for instance $V^{p,q}$ has the property that it is preserved
by $\nabla$ \ie if it is a flat subbundle, it is itself a complex
variation.

One says that $V$ is {\it polarized} by a bilinear form $b$ if $b$ is
preserved by $\nabla$ and the two Riemann bilinear equations are
verified:
$$
\eqalign{
b(u,v)&=0 \quad u\in H^{p,q}, v\in H^{r,s},\, (u,v)\not= (s,r)\cr
h_C(u,v)& =(\ii)^w b(Cu,\bar v) \,\hbox{\rm is a positive definite
hermitian metric.}\cr 
}
$$
Here $C$ is the Weil-operator which equals multiplication with
$\ii^{p-q}$ on $H^{p,q}$. Instead of $b$ one can also consider the
hermitian form 
$$
h(u,v)= \bigoplus_p (-1)^p h_C(u,v)= \ii^wb(u,\bar
w),
$$
preserved by $\nabla$. As before, we shall transport $b$, $h$ and $h_C$
to $V_\hdg$ using the $C^\infty$-isomorphism $V\mapright{\sim} V_\hdg$.

If we decompose 
$$
\nabla=\sigma^-+D+\sigma^+
$$
according to the transversality condition (*), the operator $D$ is a connection on
$V_\hdg$ and the operators $\sigma^+$ and $\sigma^-$ are
$\AA^0_M$-linear. Moreover, the operator $\sigma^+$ is the $h_C$-conjugate of
$\sigma^-$ so that one may write
$$
\eqalign{
\sigma:=& \sigma^-\cr
\sigma^*=&\sigma^+\cr 
}
$$
The $\dbar$-operator for the holomorphic structure on $V_\hdg$ is
given by the $(0,1)$-part $D^{0,1}$ of $D$ and $D$ preserves the
metric $h_C$.  So $D$ is the {\it Chern connection}, \ie the unique
metric connection on the Hodge bundle whose $(0,1)$-part is the
$\dbar$-operator.

Decomposing the equation
$\nabla\comp\nabla=0$ into types yields various equalities. The first
$$
0=\dbar(\sigma):=\dbar\comp\sigma+\sigma\comp\dbar
$$
says that $\sigma$ is a {\it holomorphic} endomorphism. The second
$$
\sigma\comp\sigma=0
$$
implies that the pair $(V_\hdg,\sigma)$ is a so-called {\it Higgs
bundle}. By definition this is a pair $(E,\theta)$ consisting of a
holomorphic bundle $E$ with a holomorphic map $\theta: E\to
\Omega^1(E)$ satisfying $\theta\wedge \theta=0$ in
$\Omega^2\End(E)$. The endomorphism $\theta$ is also called the {\it
Higgs field}. Recall also that a hermitian metric $k$ on a Higgs
bundle is called {\it harmonic} if its Chern connection $D_k$ combines 
with $\theta$ and its $k$-conjugate $\theta^*$ to give a flat
connection $\theta+D_k+\theta^*$. In our case, the Hodge metric $h_C$ is indeed
harmonic. 

Summarizing the preceding discussion, from a polarized complex
variation of Hodge bundles over $M$ we have constructed a Higgs bundle
equipped with a harmonic metric.  If $M$ is compact, such bundles can
be shown to be semi-stable.  Let me recall that a harmonic Higgs
bundle $(V,h)$ is {\it stable} resp.  {\it semi-stable} if for any
{\it proper} Higgs {\it subsheaf} $W\subset V$ \ie a coherent
submodule preserved by the Higgs field, one has an inequality of
slopes 
$$
\mu (W) < \mu(V),\,\hbox{\rm resp.}\, \mu(W)\le \mu(V).
$$
The {\it slope} for vector bundle $E$ on a K{\"a}hler manifold
$(M,\omega)$ is 
$$
\mu(E)=\deg(E)/\rank(E),
$$
with the degree of $E$ is defined using the K{\"a}hler metric:
$$
\deg(E)= c_1(E)\cdot [\omega]^{m-1},\qquad m=\dim
M.
$$
Note that $\deg(E)=\deg(\bigwedge^rE)$, $r=\rank E$.  For a torsion
free coherent sheaf $E$, one has to replace $E$ by the line bundle
which is the double dual of $\bigwedge^rE$.

Semi-stability of harmonic Higgs bundles can be proved as a consequence 
of the curvature formula, which is a straightforward consequence of
the transversality relation: 

\proclaim{Lemma} The curvature of the Chern connection $D$ on $V_{\rm
Hodge}$ is given by: 
$$ 
F_D=-[\sigma,\sigma^*].  
$$
\endproclaim

Below (see Lemma 1.4) I give a proof of a refined version of the semi-stability 
property.  Let me observe however that semi-stability is part of the
complete characterization of harmonic Higgs bundles as found by Simpson 
(see [Simpson92]):

\nnproclaim{Theorem} A Higgs bundle over a {\bf compact} K{\"a}hler
manifold $M$ with a harmonic metric is the direct sum of stable Higgs
bundles with the same slope.  Any local system which is the direct sum
of irreducible local systems admits the structure of a Higgs bundle
with a harmonic metric.  The category of Higgs bundles admitting a
harmonic metric is equivalent to the category of semi-simple local
systems, \ie those that are direct sums of irreducible local systems.
\endnnproclaim

Observe that a Higgs bundle has zero Chern classes since it carries a
flat connection and hence semi-stability means that any Higgs
subbundle has non-positive first Chern class.  I need a refinement of this
in terms of the {\it first Chern forms} 
$$
\gamma_1(E,h)={\ii \over 2\pi} \cdot \hbox{\rm (Trace of the curvature of the Chern 
  connection)}.
$$ 
So I need to estimate the curvature (with respect to the Hodge
metric) of a {\it graded Higgs subbundle} $W=\bigoplus_p W^{p,w-p}$,
\ie $W^{p,w-p}\subset V^{p,w-p}$ and $\sigma W^{p,w-p}\subset
\Omega^1_M(W^{p-1,w+p+1})$. The curvature estimate we are after reads:

\proclaim{Lemma} For a subsystem $W\subset V$ we have
$$
\ii \trace [\sigma|W_\hdg,(\sigma|W_\hdg)^*] \ge 0
$$
with equality everywhere if and only if the orthogonal complement $W_\hdg\perp$ with
respect to $h_C$ is preserved by $\sigma$ so that we have a direct sum 
decomposition of Higgs bundles
$$
W_\hdg\oplus W_\hdg\perp,
$$
or, equivalently, of complex systems of Hodge bundles.
\endproclaim
\proof If we split $V=W\oplus W^\bot$ into a $C^\infty$-orthogonal
sum with respect to the Hodge metric and write $\sigma$ into
block-form, the fact that $\sigma$ preserves $W$ means that this
block-form takes the shape
$$
\sigma=\pmatrix{ S & T \cr
 0 & S'}.
$$
Then $\trace [\sigma|W,(\sigma|W)^*]= \trace (S\comp S^*-S\comp
S^*+T\comp T^*)= \trace T\comp T^*$.  If we multiply this with $\ii$
this is  a positive definite $(1,1)$-form. The trace is zero if and
only if $T=0$ which means that $\sigma$ preserves also $W\perp$.

The last clause is a consequence of  the following discussion.
\endproof

Let me compare the flat structure $(V,\nabla)$ and the holomorphic
structure $(V_\hdg,\sigma)$ using the {\it principle of
pluri-subharmonicity} : on a compact complex manifold there are only
constant pluri-subharmonic functions. The result is:

\proclaim{Lemma}  Suppose $M$ is a smooth projective variety. A
holomorphic section of $V_\hdg$ satisfies $\sigma(s)=0$ if and only if
$\nabla(s)=0$. In other words the flat sections are precisely the holomorphic
sections of  the Hodge bundle killed by the Higgs field. 
\endproclaim
\proof 
The well known Bochner-type formula (see [Schmid, \S 7])
$$
\eqalign{
\del\dbar h_C(s,s)=& h_C(D_hs,D_hs)-h_C(F_hs,s) \ge
h_C([\sigma(s),\sigma^*(s)],s)\cr
=& -h_C(\sigma(s),\sigma(s))+h_C(\sigma^*(s),\sigma^*(s))  
}$$
shows that when $\sigma(s)=0$, the function $h_C(s,s)$ is
pluri-subharmonic and so constant. Hence $D_hs=0=\sigma^*(s)$, implying
$\nabla(s)=0$. Conversely, a flat holomorphic section satisfies
$\sigma(s)=0$ by type considerations.
\endproof

\proclaim{Corollary} A complex subbundle $W\subset V$ is a {\it
  subsystem}, \ie $\nabla W\subset W\otimes\Omega^1_X$ if and only if
$W_\hdg\subset V_\hdg$ is a Higgs subbundle. 
\endproclaim

Noticing the minus sign in the curvature of the Hodge metric and
recalling that the curvature decreases on subbundles with equality if 
and only if the quotient bundle is holomorphic, Lemma 1.2 implies:

\proclaim{Corollary} The first Chern form of a graded Higgs subbundle 
$W$ of $V$ (with respect to the Hodge metric) is negative
semi-definite and it is zero everywhere if and only if $W\perp$ is a
(holomorphic) graded Higgs subbundle as well.  In this case the
complex variation splits as {\it complex polarized variations of Hodge
structures} $V=W\oplus W\perp$.  \endproclaim

\pproclaim{Remark} These point-wise estimates can be 
integrated over any compact K{\"a}hler variety $M$ showing that the slope
of a Higgs subbundle is non-positive thereby proving semi-stability.
\endpproclaim

One also needs to know what happens for a morphism $f: V_1\to V_2 $ of
graded Higgs bundles. In general neither the kernel nor the image
are bundles, although they are preserved by the Higgs fields. Since
the kernel  of $f$ is torsion free, its degree is well-defined. For the 
image this is not necessarily the case. We have to replace it by its
saturation.  Recall that a subsheaf $F\subset V$ of a locally free (or
torsion free) sheaf $V$ is saturated if $V/F$ is torsion free. Any
subsheaf $F\subset V$ is a subsheaf of a unique saturated subsheaf
$F^{\rm sat}$, its {\it saturation}  which is defined as the inverse image
under $V\to V/F$ of the torsion of the target. 
For Higgs bundles on {\it curves} the natural morphisms
$$
\sigma_p :V^{p,w-p}\to V^{p-1,w+1-p}\otimes\Omega^1_M
$$
then can be used to define image bundles
$$
\Im(\sigma_p)= \left(\sigma_pV^{p,w-p}\right)^{\rm sat}\otimes
  (\Omega^1_M)^{\ve}.
$$

\section{A bound for variations over a compact curve}

\noindent In this section I investigate polarized complex systems of
Hodge bundles over a compact curve.  The bundles $V^{p,q}$ can be
given holomorphic structures via the $C^\infty$-isomorphism
$V^{p,q}\mapright{\sim} F^p/F^{p+1}\subset V_\hdg$ and in the sequel
we tacitly make this identification.

\proclaim{Theorem} Let $V=\bigoplus_p V^{p,w-p}$ be a weight $w$
polarized complex system of Hodge over a compact curve $C$ of genus $p$. 
Put
$$
\chi(C)=-\deg(\Omega^1_C)=2-2p.
$$
Let $\sigma=\bigoplus \sigma_p$ be the associated Higgs field and let
$\sigma^k:V^{p,w-p}\to V^{p-k,w-p+k}$ the composition
$\sigma_{p-k+1}\comp\cdots\comp \sigma_p$. Then we have the inequality
$$
\deg V^{p,w-p} \le -\chi(C)\sum_{r\ge 1} {r\over 2} \left(\rank
  \sigma^{r-1}-\rank \sigma^r\right)
$$
Equality holds if and only if for some $k\ge 0$ the maps
$\sigma_p,\dots,\sigma_{p-k}$ are all isomorphisms and $\sigma_{p-k-1}$
is zero. In this case $V^{p,w-p}\oplus\cdots\oplus V^{p-k-1,w-p+k+1}$
is a  complex subvariation. 
\endproclaim
 
\noindent{\bf Proof (Deligne)}:
Consider the following graded Higgs subbundles $W_r\subset V$, $r=1,2,\dots$ which
only are non-zero in degrees $p,\dots, p-r+1$ and which  --- using the notation
of the image under $\sigma$ as defined in the previous section ---
are given by
$$
\eqalign{
 W_r^{p,w-p}:=\ker\sigma^r\subset V^{p,w-p}, \quad
 &W_r^{p-1,w-p+1}:=\Im(\sigma\ker\sigma^r)\subset V^{p-1,w-p+1},\dots
 \cr \, \dots &\quad
 W_{r}^{p-r+1,w-p+r-1}:=\Im(\sigma^{r-1}\ker\sigma^r)\subset
 V^{p-r+1,w-p+r-1}.  
 }
$$
Apply the stability property to these bundles. For simplicity we put
$$
\eqalign{
d_r& = \deg (\ker\sigma^r /\ker \sigma^{r-1})\cr
l_r& = \dim (\ker\sigma^r /\ker \sigma^{r-1})= \rank
\sigma^{r-1}-\rank \sigma^r
}
$$
Now compute the degree of each  of these $r$ bundles 
$W_{r}^{p-k,w-p+k}$ using the isomorphism
$$
\Ker \sigma^r/\ker \sigma^k \mapright{\sim}  \sigma^k(\Ker
\sigma^r)^{\rm sat}\subset V^{p-k,w-p+k}\otimes \left(\Omega^1_C\right)^{\otimes
k}.
$$
The bundle on the left has rank $l_r+\cdots +l_{k+1}$ and degree
$d_r+\cdots+d_{k+1}$.  It follows that $\deg
W_{r}^{p-k,w-p+k}=\deg\left(\Im(\sigma^k\Ker \sigma^r)\right)=
\deg\left((\Ker\sigma^r/\ker \sigma^k) \otimes \left(\Omega^1_C\right)^{-\otimes
k}\right) = d_r+\cdots+d_{k+1}+(l_r+\cdots +l_{k+1})k\chi$ and adding these for
$k=0,\dots,r-1$ one finds the total degree 
$$
\deg W_r= \left(  \sum_{p=1}^rpd_p+l_p{p(p-1)\over 2}
\cdot \chi(C) \right)  \le 0.
$$
Now take a weighted sum of these inequalities in order to make appear
$\deg V^p=\sum d_p$. Indeed
$$
 0\ge \sum _{k=1}^\infty
\left( {1\over k}-{1\over k+1}\right) \cdot\left(\sum_{p=1}^k
  pd_p+l_p{p(p-1)\over 2}\cdot \chi(C) \right) =\sum_{p\ge 
  1}d_p + \sum_{p\ge 1}  l_p  {p-1 \over 2}\cdot \chi(C) 
$$
You have equality if and only if all of the equalities you started out
with are equalities which means that the $W_r$ are direct factors as
graded Higgs subbundle of $V$. This translates into the $\sigma_r$
being either zero or an isomorphism. To show this, consider the
bundle $W_1$, \ie $\Ker\sigma_p\subset V^{p,w-p}$. If this is direct
factor of $V$ as a graded Higgs  bundle, the orthogonal complement has
the structure of a graded Higgs bundle. Since the latter, if not zero,
must have type different from $(p,w-p)$ it follows that $\sigma_p$ is
either $0$ or has maximal rank.  If it is the zero map, the kernel $V^{p,w-p}$ is
itself a subvariation.  If it has maximal rank, its kernel is zero and
hence has degree $0$ and we repeat the argument with $\sigma_{p-1}$. 
If this map is zero $V^{p,w-p}\oplus V^{p-1,w-p+1}$ is a subvariation
and if $\sigma_{p-1}$ is an isomorphism we continue the procedure.  In
this way we produce a chain $\sigma_p,\dots,\sigma_{p-k}$ of
isomorphisms such that the next map $\sigma_{p-k-1}$ is zero as
asserted.  \endproof

Next, we remark that these bounds applied to the {\it dual} variation
of Hodge structure gives bounds in the other direction. We shall make
these explicit for real variations, which are self-dual in the obvious 
sense. In fact, the reality constraint $H^{p,q}=\bar H^{q,p}$ implies
$\deg H^{p,q}=-\deg H^{q,p}$. In the real case, again one has equality
only for a consecutive  string $\sigma_p,\dots,\sigma_{p-k}$ of
isomorphisms such that the next one is zero. Note however that the complex
subvariation splitting the real variation itself need not be real. 
The reality constraint itself imposes however some additional constraints.

\proclaim{Corollary}
Let $V$ be a variation of polarized real weight $w$ Hodge structures over a
compact curve $C$ with Higgs field $\sigma=\bigoplus\sigma_p$. Put
$\bar\sigma_p=\sigma_{w-p}$ and
$\sigma^k=\sigma_{p-k-1}\comp\cdots\comp \sigma_p:V^{p,w-p}\to
V^{p-k,w-p+k}$ and $\bar\sigma^k=\bar\sigma_{p-k-1}\comp\cdots\comp
\bar\sigma_p:V^{w-p,p}\to V^{w-p+k,p-k}$.
Then we have the bounds
$$
\eqalign{\chi(C)\sum_{r\ge 1} {r\over 2} \left(\rank
 \bar\sigma^{r-1}-\rank \bar\sigma^r\right)&\le\deg V^{p,w-p}\le \cr
&-\chi(C)\sum_{r\ge 1} {r\over 2} \left(\rank\sigma^{r-1}-\rank \sigma^r\right).
}$$
In particular, if $V^{p,w-p}\not= 0$ precisely in the interval
$p=0,\dots, w$, the Hodge bundle $V^{w,0}$ satisfies the bounds
$$
0 \le \deg V^{w,0} \le -\chi(C)\sum_{r=1}^w {r\over 2}
\left(\rank\sigma^{r-1}-\rank \sigma^r\right)
$$
and equality on the left holds if and only if $V^{w,0}$ is a flat
subbundle, while equality on the right holds if and only if {\rm all}
the maps $\sigma_k$ are isomorphisms. 
\endproclaim

Let me translate this Corollary in terms of period maps. The period 
domains involved are the Griffiths domain $D$ resp. $D'$ parametrizing real polarized 
Hodge structures of weight $w$ of fixed given type, resp.  the partial Hodge flags
$V^{w,0}\subset V$ of the same type.  Giving a polarized variation of
Hodge structures of this type over the curve $C$ is the same as giving
its period map
$$
p: C \to D/\Gamma,
$$
where $\Gamma$ is the monodromy group. Let me refer to [Griffiths] for 
the necessary background. Likewise, the partial Hodge flag is given 
by the partial period map
$$
p':C\to D'/\Gamma' 
$$ defined as the composition of $p$ and the forgetful map
$q:D/\Gamma\to D'/\Gamma'$.  The subbundle $V^{w,0}=F^w$ being flat
means exactly that $p'$ is constant. At the other end of the spectrum, 
using that the maps $\sigma_k$ measure the derivative of the period
map, the latter is everywhere an immersion if and only if all
$\sigma_{k}$ are isomorphisms.  Summarizing: 

\proclaim{Corollary} 
\item{i)} One has $\deg W^{w,0}=0$ if and only if the image of the full period map
lands in one of the fibers of $q$. 
\item{ii)} The equality $\deg V^{w,0}= -\chi(C)\sum_{r=1}^w {r\over 2}
\left(\rank\sigma^{r-1}-\rank \sigma^r\right)$ holds if and only if 
the perod map is everywhere an immersion.

\endproclaim

\pproclaim{Remark} 
Put
$$
\eqalign{
v^{p,q}:=&\dim V^{p,q}\cr
v_0^{p,q}:=&\dim \ker (\sigma_p:V^{p,q}\to V^{p-1,q+1}).\cr
}
$$
Using the estimate
$$
\rank \sigma_{w-r}\le \rank\sigma^{r+1}
$$
and a telescoping argument, the  above  inequality implies the (a
priori weaker) bound 
$$
\deg E^{w,0}\le -\chi(C) \cdot \sum_p (v^{p,w-p}-v_o^{p,w-p})
$$
and this is due to [JostZuo] (they combine this with the  obvious
symmetry among the Hodge numbers and, similarly,  the numbers $v^{p,q}_0$).
\endpproclaim

\section{The case of a smooth non-compact curve}

We let $C$ be a smooth compactification and we let $S=C\setminus
C_0$.  Assume now that $V_0$ is a complex variation of Hodge structure
with the property that the local monodromy operators $\gamma_s$ around
the points $s\in S$ are quasi-unipotent.  This is true if the complex
variation is defined over the integers, e.g. if the variation comes
from geometry.  This is the monodromy theorem, whose proof can be
found in [Schmid]. 

The assumption imply that $V_0$ admits a vector bundle extension $V$ 
over $C$ and all the Hodge bundles $F^pV_0$ extend to holomorphic
bundles $F^pV\subset V$ on $C$. The Higgs fields are now maps
$$
\sigma: V\to V\otimes \Omega_C^1(\log S).
$$
Asymptotic analysis of the Hodge metric $h_p$ on every Hodge bundle
$F^p$ around the punctures shows that its Chern form is integrable and
that there are non-negative rational numbers $\alpha_s^p$, the {\it
residues} of the operators $\gamma_p$ such that
$$
\deg(F^p)=\int_C \gamma_1(F^p,h_p)+\sum_{s\in S}\alpha_s^p.
$$
For details of this I refer to [Peters84]. The residues are always
zero in the unipotent case and conversely, if the Hodge bundles are
{\it flat} the vanishing of all the residues implies that the local
monodromy operators are all unipotent. In general, the residue is
given by
$$
\alpha_s^p=\sum_{\alpha\in [0,1)} \alpha \dim (F^p\cap V_\alpha)
\leqno{\rm (*)}
$$
where $V_\alpha$ is the subspace of $V$ on which the local monodromy
operator $T_s$ acts with eigenvalue $\exp(2\pi i\alpha)$.

The estimates of the previous section go through for the first summand 
and also for its successive quotients, the Hodge bundles. So we find

\proclaim{Theorem} Let $V$ be  a complex variation over a curve
$C_0=C\setminus S$ with quasi-unipotent local monodromy operators. Put
$\genus(C)=g$. Let $\sigma=\bigoplus \sigma_p$ be the associated Higgs
field and let $\sigma^k:V^{p,w-p}\to V^{p-k,w-p+k}$ the composition
$\sigma_{p-k+1}\comp\cdots\comp \sigma_p$.
Then the degree of the Hodge bundles are bounded  by 
$$
\deg V^{p,w-p} \le (2g-2+\#S)\cdot \sum_r {r\over 2} \left(\rank
  \sigma^{r-1}-\rank \sigma^r\right) +\sum_{s\in S}
(\alpha^{p}_s-\alpha^{p+1}_s), 
$$
with $\alpha_s$ given by {\rm (*)}.
\endproclaim

If the local monodromy-operators are unipotent, the last summand 
vanishes and one has a description for those variations where the 
bound is attained:

\proclaim{Addition} In the situation of the preceding theorem, assume 
that the local monodromy operators are unipotent.  Equality is
attained if and only if for some $k\ge 0$ the maps
$\sigma_p,\dots,\sigma_{p-k}$ are all isomorphisms and
$\sigma_{p-k-1}$ is zero.  In this case $V^{p,w-p}_{0}\oplus\cdots\oplus
V^{p-k-1,w-p+k+1}_{0}$ is a complex subvariation of $V_{0}$.
\endproclaim
The proof is the same as before; it depends on the validity of Lemma 
1.3 in the case of non-compact curves. The  principle of
pluri-subharmonicity needed here can be found in [Schmid], proof
of Theorem (7.22). It really depends on his $\SL(2)$-orbit theorem and
so it is considerably less trivial.

\pproclaim{Remark}
For real variations there are bounds as in Corollary 2.2:
$$
0\le \sum _{s\in S} \alpha_s^p \le \deg V^{w,0} \le (2p-2+\#S)\cdot
\sum_{r=1}^w {r\over 2} \left(\rank\sigma^{r-1}-\rank
\sigma^r\right)+\sum _{s\in S} \alpha_s^p 
$$
and if for instance $\deg V^{w,0}=0$, all the local monodromy operators 
are unipotent and $V^{w,0}$ is a flat subbundle of $V$. See also [Peters84].

\endpproclaim

\section{Boundedness revisited}

Here we give the proof that a bound on the degrees of all of the Hodge 
bundles imply boundedness. To be precise, we need both upper and
lowerbounds so that only finitely many degrees are possible. 
We repeat here that lower bounds are obtained from the upper bounds of
the dual variation.

First we need to see that only those variations that occur on a fixed
local systems need to be considered. For this we invoke Deligne's
boundedness result from [Deligne]: 
\nnproclaim{Theorem} Given a smooth algebraic variety $S$ there are
only finitely many   classes of representations of $\pi_1(S)$ on a
given rational vector  space $V$ such that the resulting flat vector
bundle underlies a  polarizable complex variation of Hodge  structure of
a given weight. 
\endnnproclaim

I would like to point out that the proof of this theorem is not too 
difficult, except for one point: it uses the (by now standard fact)
that a polarizable complex variation of Hodge  structure is
semi-simple; in any case it is much easier than Simpson's methods used
to show boundedness. 

We complete the argument by showing:

\proclaim{Lemma} Let there be given a flat vector bundle on a curve $C$ and a 
finite set  of integers $\{d_{p}\}$, $p\in I$. Complex variations 
$V=\bigoplus_{p\in I} V^{p,w-p}$ with $\deg V^{p,w-p}=d_{p}$ are 
parametrized by an open set of an algebraic variety.   
\endproclaim
\proof Consider the bundle $\FF$ over $C$ whose fiber over $c$ consists of all 
flags in $V_c$, the fiber of the given complex variation $V$ at $c$,
which are of the same type as the Hodge flags. Then the total space
$\FF$  is a projective variety. There is a subvariety $\GG$ of $\FF$
formed by flags   satisfying the first bilinear relation. Over $\FF$
we have the  tautological flag  of vectorbundles $F^{p}$. A complex variation 
underlying $V$ is the same thing as a section $s$ of $\GG$ which in each 
fiber lands in the set of flags satisfying the second period 
relation. The degree of the $p$-th Hodge bundle is the degree of 
$s^{*}F^{p}$, together making up the {\it flag degree} $\{d_{p}\}$, 
$p\in I$. Since these are all fixed, the desired collection of 
complex variations forms an open set in the variety $S$ of sections of the 
partial flag bundle $\GG$ over $C$ of fixed flag degree. This 
variety $S$  can be viewed as  the Hilbert scheme of curves in 
$\GG$ of fixed multidegree with respect to the embeddings defined by 
the flags.
\endproof

\bigskip

\centerline{\sectionfont References } 
\def\refitem[#1]#2:#3#4#5(#6)#7%
{\hangindent=6.85em\hangafter=1
\noindent\rlap{\hbox{[{\bf #1}]}}\kern6.5em{\rm #2: }%
     {\rm #3, }{\rm #4 }{\bf #5 }({\rm #6}) {\rm #7}\par}

\refitem[Arakelov] Arakelov, A.:{Families of algebraic curves with fixed
  degeneracies}{Izv. Ak. Nauk. S.S.S.R,
  ser. Math.}{35}(1971){1277--1302.} 

\refitem[Deligne] Deligne, P.:{Un th{\'e}or{\`e}me de finitude pour la
monodromie, in {\sl Discrete groups in geometry and analysis,
Papers in honor of G. D. Mostow's 60th
Birthday}}{Prog. Math.}{67}(1987){1--19.}   

\refitem[Eyss] Eyssidieux, Ph.:{La caract{\'e}ristique d'Euler du
complexe de Gauss-Manin}{Journ. f. reine und angew. Math.}{490}(1997){155--212.}

\refitem[Faltings] Faltings, G.:{Arakelov's theorem for abelian
  varieties}{Invent. Math.}{73}(1983){337--348.}   

\refitem[Griffiths] Griffiths, P. A.:{\it Topics in transcendental
algebraic geometry}{Princeton Univ.  Press}{}(1984){.}

\refitem[JostZuo] Jost, J., K. Zuo:{Arakelov type inequalities for Hodge
  bundles over algebraic varieties, Part I: Hodge bundles over
  algebraic curves with unipotent monodromies around
  singularities}{Preprint}{}(1999){}

\refitem[Peters84] Peters, C.A.M. :{A criterion for flatness of Hodge
  bundles over curves and geometric
  applications}{Math. Ann.}{268}(1984){1--19.}

\refitem[Peters86] Peters, C.A.M.:{On Arakelov's finiteness theorem
  for higher dimensional
  varieties}{Rend. Sem. Mat. Univ. Politec. Torino}{}(1986){43--50.}

\refitem[Peters90] Peters, C.A.M. :{ Rigidity for variations of Hodge
structure and Arakelov-type finiteness theorems}{Comp. Math}{75}(1990){113--126.}

\refitem[Saito] Saito, M.-H.:{Classification of non-rigid families
  of abelian varieties}{Tohoku Math. J}{45}(1993){159--189.}  

\refitem[SaitoZucker] Saito, M.-H., S. Zucker:{Classification of
  non-rigid families of K3-surfaces and a finiteness theorem of
  Arakelov type}{Math. Ann. 289}(1991){1--31.}

\refitem[Schmid] Schmid, W.:{Variation of Hodge structure: the
  singularities of the period
  mapping}{Invent. Math.}{22}(1973){211--319.}

\refitem[Simpson90] Simpson, C.:{Harmonic bundles on non-compact
  curves}{J. Am. Math. Soc.}{3}(1990){713--770.} 

\refitem[Simpson92] Simpson, C.:{Higgs bundles and local
  systems}{Publ. Math. I.H.E.S.}{75}(1992){5--95.}

\refitem[Simpson94] Simpson, C.:{Moduli of representations of the
  fundamental group of a smooth variety}{Publ. Math. I.H.E.S.}{79}(1994){47--129.}

\refitem[Simpson95] Simpson, C.:{Moduli of representations of the
  fundamental group of a smooth projective
  variety. II}{Publ. Math. I.H.E.S.}{80}(1995){5-79.}

\bigskip
\begingroup
\smallfonts

\endgroup

\end